\documentclass[12pt]{article}
\usepackage{amscd,amsthm}
\usepackage{latexsym}
\usepackage{amsmath}
\usepackage{amssymb}

\def\bbR{{\Bbb R}}
\def\bbZ{{\Bbb Z}}

\def\f{{\cal F}}

\def\x{{\cal X}}

\def\const{{\rm const}}

\def\tr{{\rm Tr}}
\def\id{{\rm Id}}
\def\var{{\rm Var}}

\begin{document}
\title{Gaussian Limit for Determinantal Random Point Fields
\footnote{ AMS 2000 subject classification: 60G55 (primary), 60F05 (secondary);
keywords and phrases: determinantal random point fields, Central Limit Theorem,
self-similar random processes}}
\author{Alexander Soshnikov\\
University of
California, Davis 
\\Department of Mathematics\\Davis, CA  95616, USA\\
soshniko@math.ucdavis.edu}
\date{}
\maketitle
\begin{abstract}
We prove that, under fairly general conditions, a properly
rescaled determinantal random point field converges to
a generalized Gaussian random process.
\end{abstract}

\section{Introduction and Formulation of Results}

Let $E$ be a locally compact Hausdorff space satisfying the
second axiom of countability, $B-\sigma$-algebra of Borel
subsets and $\mu \ \ \ $a $\sigma$-finite measure on $(E,B)$, such
that $\mu (K)<\infty$ for any compact $K\subset E$.  We
denote by $X$ the space of locally finite configurations
of particles in $E$:  $X=\{\xi =(x_i)^\infty_{i=-\infty}:
x_i\in E\ \forall i$, and for any compact $K\subset E$
$\#_K(\xi ):=\#(x_i:x_i\in K)<+\infty\}$.  A $\sigma$-algebra
$\f$ of measurable subsets of $X$ is generated by the
cylinder sets $C^B_n=\{\xi\in X:\#_B(\xi )=n\}$, where $B$
is a Borel set with a compact closure and $n\in\bbZ^1_+=
\{0,1,2,\dots\}$.  Let $P$ be a probability measure on
$(X,\f )$.  A triple $(X,\f ,P)$ is called a random point
field (process) (see [DVJ], [Le1-Le3]).  In this paper we
will be interested in a special class of random point fields
called determinantal random point fields.  It should be
noted that most, if not all the important examples
of determinantal point fields arise when $E=\coprod^k_{i=1}
E_i \ $ (here we use the notation $ \coprod $  for the disjoint union),
$\ E_i\cong\bbR^d$ or $\bbZ^d$ and $\mu$ is either the
Lebegue or the counting measure.  We will, however, develop our
results in the general setting ( our arguments
will not require significant changes).

Let $dx_i,\ i=1,\dots ,n$ be disjoint infinitesimally small subsets
around the $x_i$'s.  Suppose that a probability to find a particle in
each $dx_i$ (with no restrictions outside of $\coprod^n_{i=1}
dx_i$) is proportional to $\prod^n_{i=1}\mu (dx_i)$, i.e.
\begin{equation}
P(\#(dx_i)=1,\ i=1,\dots ,n)=\rho_n(x_1,\dots ,x_n)\mu (dx_1)\dots
\mu (dx_n)
\end{equation}
The function $\rho_n(x_1,\dots ,x_n)$ is then called the $n$-point
correlation function.  The equivalent definition is given by the
equalities
$$E\prod^m_{i=1}\frac{(\#_{B_i})!}{(\#_{B_i}-n_i)!}=\int_{B^{n_1}_1
\times\dots\times B^{n_m}_m}
\rho_n(x_1,\dots ,x_n) d\mu (x_1)\dots 
d\mu (x_n)$$ 
where $B_1,\dots ,B_m$ are disjoint Borel sets with
compact closures, $m\geq 1,\ n_i\geq 1,\ i=1,\dots ,m, n_1+\dots
+n_m=n$.

A random point field is called determinantal if
\begin{equation}
\rho_n(x_1,\dots ,x_n)=\det\bigl (K(x_i,x_j)\bigr )_{1\leq 
i,j\leq n},
\end{equation}
where $K(x,y)$ is a kernel of an integral operator $K: L^2(E,d\mu )
\rightarrow L^2(E,d\mu )$ and $K(x,y)$ satisfies some natural 
regularity conditions discussed below. Such a kernel
$K(x,y)$ is  called a correlation kernel.

It follows from 2) and the non-negativity of the $n$-point
correlation functions that $K$ must have non-negative minors, and
in particular if $K$ is Hermitian it must be a non-negative operator.
In this paper we shall always restrict ourselves to the Hermitian
case.

Determinantal (also known as fermion) random point fields were introduced
by Macchi in the early seventies (see [Ma1], [Ma2], [DVJ]).  A recent
survey of the subject with applications to random matrix theory,
statistical mechanics, quantum mechanics, probability theory and
representation theory is given in [So1].  Diaconis and Evans in [DE1]
introduced a generalization of determinantial random point processes,
called immanantal point processes.

Let $K$  be a Hermitian,
locally trace class, integral operator on $L^2(E,d\mu )$.  Suppose
that we can choose a kernel $K(x,y)$ in such a way that for any
Borel set $B$ with compact closure
\begin{equation}
\tr (K\x_B)=\int_BK(x,x)dx,
\end{equation}
where $\x_B$ denotes the multiplication operator by the indicator of
$B\ (\equiv$projector on the subspace of the functions supported in
$B$).

Since it is always true that
\begin{gather*}
\tr (K\x_{B_1}\dots K\x_{B_n})=
\\ \int_{B_1\times\dots\times B_n}
K(x_1,x_2)K(x_2,x_3)\dots K(x_n,x_1)
d\mu (x_1)\dots d\mu (x_n)
\tag{3'}
\end{gather*}
for $n>1$ and Borel sets $B_1,\dots ,B_n$ with compact closure, (3)
implies that (3$'$) holds for all $n$.

(3) can always be achieved for $E=\bbR^d$ (see e.g. [So1], Lemmas
1,2).  From now on we will assume that both (2) and (3) are satisfied.

The main goal of our paper is to study the behavior of linear
statistics
$$S_f(\xi )=\sum_i f(x_i),\ \xi =(x_i),$$
for sufficiently ``nice" test functions in a scaling limit.
The moments of $S_f$ can be calculated from (2).  For instance,
\begin{equation}
ES_f=\int f(x) K(x,x) d\mu (x),
\end{equation}
\begin{equation}
{\rm Var}\ S_f=\int f^2(x)K(x,x)d\mu (x)-\int f(x)f(y)\vert K(x,y)
\vert^2d\mu (x)d\mu (y).
\end{equation}
Taking $E=\bbR^1$ and $K(x,y)=\tfrac{\sin\pi (x-y)}{\pi (x-y)}$, a
so-called sine kernel, we obtain a random point field well known
in the theory of random matrices.  It can be viewed as a limit
$n\rightarrow\infty$ of the distribution of the appropriately scaled
eigenvalues of $n\times n$ random Hermitian matrices with Gaussian
entries (see e.g. [D], chapter 5).  It was proven by Spohn in [Sp]
(see also [So2]), that if $K$ is the sine kernel and a test
function $f$ is sufficiently smooth and fast decaying at infinity, 
then $\sum^\infty_{i=-\infty}f(\tfrac{x_i}{L})
-L\int^\infty_{-\infty}f(x)dx$ converges in distribution
to the normal law $N(0,\int^\infty_{-\infty}\vert\hat f
(k)\vert^2\cdot\vert k\vert dk)$, where $\hat f$ is the
Fourier transform of $f$, $\hat f(k)=\int^\infty_{-\infty}
f(x)e^{-2\pi ikx}dx$.  In other words we can say that the
random signed measure 
$$\sum^\infty_{i=-\infty}\delta (x-\tfrac{x_i}{L})-Ldx$$
converges as $L\rightarrow\infty$ to the generalized
self-similar Gaussian random process with the spectral
density $\vert k\vert$ (see e.g., [Dob1], [Dob2] \S 3, 
[S], and, for the introduction to the theory of generalized
random processes, [GV]).  The fact that the variance of the
linear statistics $\sum^\infty_{i=-\infty}f(\tfrac{x_i}
{L})$ does not grow to infinity for Schwartz functions
is the manifestation of the strong repulsiveness of the
distribution of the eigenvalues of random matrices.  Similar
results for other ensembles of random matrices have been
obtained in [DS], [Jo1], [Jo2], [B], [BW], [So2], [W],[DE2].  The
kernels appearing in these ensembles are, in some respect,
very much like the sine kernel.  In particular, the variance
of the number of particles in an interval grows as a 
logarithm of the mathematical expectation of the number of
particles.  The following result was established by Costin
and Lebowitz for the sine kernel ([CL]):  let $f$ be an
indicator of an interval, $f=\x_I,\ I=(a,b)$, then
$$E\sum^\infty_{i=-\infty}f(x_i/L)=E(\#(x_i: aL<x_i
<bL))=L(b-a),$$ 
$${\rm Var}\ (\sum^\infty_{i=-\infty}f(
x_i/L))=\frac{1}{\pi^2}\log L+O(1),$$ 
and 
$$\frac{\#(x_i:
aL<x_i<bL)-L(b-a)}{\sqrt{\frac{1}{\pi^2}\log L}}$$ 
converges
in distribution as $L\rightarrow\infty$ to the normal law
$N(0,1)$.  The proof of the Costin-Lebowitz theorem holds, quite
remarkably, for arbitrary determinantal random point fields
with Hermitian kernel.
\medskip

\noindent{\bf Theorem} ([So3]) {\it Let $(X,\f ,P_L), \ L\geq 0,$ 
be a family of
determinantal random point fields with  Hermitian locally trace class kernels
$K_L$ and $\{I_L\}_{L\geq 0}$ be a family of Borel subsets
of $E$ with compact closure.  Then if ${\rm Var}_L\ (\#(x_i:
x_i\in I_L))\underset{L\rightarrow\infty}{\rightarrow}
\infty$, the normalized random variable $\tfrac{\#(x_i:
x_i\in I_L)-E_L\#_{I_L}}{\sqrt{{\rm Var}_L\ \#_{I_L}}}$
converges in distribution to $N(0,1)$.}
\medskip

Here and below we denote by $E_L, \ \ {\rm Var}_L$  the mathematical 
expectation and the variance with respect to $P_L$.
One can also establish a similar result for the step
functions (finite linear combinations of indicators).
\medskip

\noindent{\bf Theorem} {\it Let $(X,\f ,P_L)$ be a family of
determinantal random point fields with  Hermitian locally trace class kernels
$K_L$ and $\{I^{(1)}_L,\dots ,I^{(k)}_L\}_{L\geq 0}$ be a
family of Borel subsets of $E$, disjoint for any fixed $L$,
with compact closure.  Then if for some $\alpha_1,\dots ,
\alpha_k\in\bbR^1$, the variance of the linear statistics
$\sum^\infty_{i=-\infty}f_L(x_i)$ with $f_L(x)=\sum^k_{j=1}
\alpha_j\cdot\x_{I^{(j)}_L}(x)$,  grows to infinity in such a way
that ${\rm Var}_L (\#_{I^{(j)}_L})=O({\rm Var}_L\ (\sum^\infty_{i=
-\infty}f_L(x_i)))$ for any $1\leq j\leq k$, the Central
Limit Theorem holds:
$$\frac{\sum^k_{j=1}\alpha_j^{(L)}\cdot\#_{I^{(j)}_L}-E_L\left (
\sum^k_{j=1}\alpha_j\cdot\#_{I^{(j)}_L}\right )}{\sqrt
{{\rm Var}_L\ (\sum^k_{j=1}\alpha_j\cdot\#_{I^{(j)}_L})}}
\overset{w}{\longrightarrow}N(0,1).$$}
\medskip

\noindent{\bf Remark 1} We use standard notations $f=O(g)$
and $f=o(g)$ when $\tfrac{f}{g}$ stays bounded or $\tfrac
{f}{g}\rightarrow 0$.
\medskip

\noindent{\bf Remark 2}  The last theorem has been
explicitly stated in [So3] only in the special case of
the Airy and Bessel kernels and the kernels arising in the 
classical compact groups (see Theorems 1, 2, 4, 6), however
the key Lemmas 7 and 8 proven there allow rather
straightforward generalization to the case of an arbitrary
Hermitian kernel.  A result close to our Theorem 6 from [So3] was
also established by K. Wieand ([W]).

We recall that a Hermitian kernel $K(x,y)$ defines a
determinantal random point field if and only if the
integral operator $K$ is non-negative and bounded from
above by the identity,
\begin{equation}
0\leq K\leq\id
\end{equation}
([So1], Theorem 3).  For the translation-invariant
kernels $K(x-y)$ and $ E= \bbR^d$  or $\bbZ^d $ 
this is equivalent to $0\leq \hat K(t)\leq
1$, where
\begin{equation}
K(x)=\int e^{2\pi i(x\cdot t)}\hat K(t)dt
\end{equation}
The sine kernel $K(x-y)=\tfrac{\sin\pi (x-y)}{\pi (x-y)}$
corresponds to $\hat K(t)=\x_{\left
[-\tfrac{1}{2},\tfrac{1}{2}\right ]}
(t)$, the indicator of $[-\tfrac{1}{2},\tfrac{1}{2}]$.
It might be worth noting and  actually is not very difficult
to see, that the logarithmic rate of the growth of Var $(\#(x_i:\vert
x_i\vert\leq L))$ for the sine kernel is the slowest among
all translation-invariant kernels corresponding to projectors,
$\hat K = \x_B $, for which $inf(B)$ and $sup(B)$ are the density points
of $B$. 
For the generic translation-invariant kernel
$K(x-y)$ ( $\hat K$ is not an indicator) Var $(\#(
x_i:\vert x_i\vert\leq L))$ is proportional to Vol $(x_i:\vert x_i
\vert\leq L)\sim E(\#(x_i:\vert x_i\vert\leq L))$ ([So1], section 3).

In our main result we prove CLT for the linear statistics
when the variance grows faster than some arbitrary
small, but fixed, power of the mathematical expectation.

\medskip

\noindent{\bf Theorem 1} {\it Let $(X, \f, P_L), \ L\geq 0 $ be a family of 
determinantal random point fields with Hermitian correlation kernels $K_L$.
Suppose that
$f_L, \ L\geq 0$ are 
bounded measurable  functions with precompact support 
(i.e.  $supp(f_L) $ has a compact 
closure for any $L \geq 0$), such that

\begin{equation}
{\rm Var}_L\ S_{f_L}\rightarrow\infty\text{ as }L\rightarrow
\infty
\end{equation}
and
\begin{equation}
sup |f_L(x)| =o ({\rm Var}_L)^{\epsilon} , \ \ 
E_LS_{\vert f\vert_L}=O\bigl (({\rm Var}_L\ S_{f_L})^\delta
\bigr ),
\end{equation}
for any $\epsilon >0$  and some  $\delta >0$.  Then the normalized linear 
statistics
$\tfrac{S_{f_L}-E_LS_{f_L}}{\sqrt{{\rm Var}_L\ S_{f_L}}}$
converges in distribution to the standard normal law
$N(0,1)$.}  
\medskip

As a very important special case of  Theorem 1 one can consider 
$f_L(x):=f(T_Lx)$, where
$\{T_L\},\ L\in\bbR^1_+$, is a one-parameter family
of measurable transformations $T_L:E\rightarrow E$ such 
that $T^{-1}_LD$ has compact closure for any compact $D$.
If for a sufficiently rich class of test
functions $f$ (e.g., continuous functions with compact support)
(8),(9) are satisfied, and
the rate of the growth of ${\rm Var}_L  \ (S_{f_L})$ is the same,
$${\rm Var}_L\ (S_{f_L})=B(f)\cdot V_L\cdot (1+o(1)),$$ 
where $B(f)$ is some functional on a space of test
functions, Theorem 1 implies that the random signed measure
$$V^{-\frac{1}{2}}_L\left (\sum^\infty_{i=-\infty}\delta
(x-T_Lx_i)-T_L (K_L(x,x)d\mu (x))\right )$$
converges as $L\rightarrow\infty$ to the generalized
Gaussian process with the correlation functional $B(f,f)=
B(f)$ ( we denote by
$T_L (K_L(x,x)d\mu (x)) $ the image of the measure
$K_L(x,x)d\mu (x) $ under $T_L$).

Let us consider a Euclidean one-particle space
$E=\bbR^d \ $, a one-parameter family of dilations
$$ T_L :\bbR^d \rightarrow \bbR^d, \ \ T_L x= x/L, $$
and a  correlation kernel
\begin{equation}
K_L(x,y)=A_L(x-y)+R_L(x,y),
\end{equation}

where
\begin{equation}
\vert R_L(x,y)\vert\leq Q(x_{abs}+y_{abs}),
\end{equation}
$x_{abs}=(\vert x_1\vert,\dots ,\vert x_d\vert ),\ Q\in
L^2(\bbR^d_+)\cap L^\infty (\bbR^d_+)$ 
It follows from (6),
(10) and (11) that $0\leq A_L\leq\id$, which implies 
$0\leq\hat A_L(k)\leq 1$, $0\leq\int_{\bbR^d}\hat A_L(k)-(\hat
A_L(k))^2dk=A_L(0)- \int_{\bbR^d} |A_L(x)|^2 dx =:\sigma_L^2$, and  $\sigma_L =0$ if and only if $\hat A_L$
is an indicator.
\medskip

\noindent{\bf Theorem 2} {\it Let the kernel $K_L$ satisfy
(10), (11) and there exist constants $const, \ \sigma >0$ and 
$\kappa_L \to \infty$ as $L \to \infty$ such that
$$ \sigma_L \to \sigma \ \ {\rm as} \  L \to \infty ,$$
$$ |A_L(0)| <const, $$
and
$$ \int_{|x|>L/\kappa_L} |A_L(x)|^2 dx \to 0.$$
Then for any real-valued
function $f\in L^1(\bbR^d)\cap L^2(\bbR^d)$ the 
normalized linear statistics
$$\frac{1}{L^{\frac{d}{2}}\sigma}\left (\sum^\infty_{
i=-\infty}f\left (\frac{x_i}{L}\right )-A_L(0)\cdot L
\int_{\bbR^d}f(x)dx\right )$$
converges in distribution to the Gaussian random variable
$N\left (0,\int_{\bbR^d}\left (f(x)\right )^2dx\right )$.}
\medskip

\noindent{\bf Remark 3} Theorem 2 says that under the stated 
conditions the random
signed measure
$$\frac{1}{L^{\frac{d}{2}}\sigma}\left (\sum^\infty_{i=-
\infty}\delta\left (x-\frac{x_i}{L}\right )-A_L(0)\cdot Ldx
\right )$$
converges to the white noise  as $L\rightarrow
\infty$ (for the definition of the white noise see e.g.
[H]).  Similar results hold in the discrete case.

Let us now restrict our attention to the 
translation-invariant kernels $K(x,y)=A(x-y)$.  We will 
use the notation
$$m(\lambda ):=\int\hat A(k)-\hat A(k)\hat A(k-\lambda )
dk.$$
Observe that $\sigma^2=m(0)$ and
\begin{equation}
{\rm Var}\ \left (\sum^\infty_{i=-\infty}f(x_i)\right )
=\int\vert\hat f(\lambda )\vert^2m(\lambda )d\lambda .
\end{equation}
In particular
\begin{equation}
{\rm Var}\ \left (\#_{[-L,L]^d}\right )={\rm Vol}\ \left
([-L,L]^d\right )\cdot\left (m(0)+o(1)\right ).
\end{equation}
It follows from (12) that the rate of the growth of the
variance of $S_{f_L}$ depends on the asymptotics of $m(
\lambda )$ near the origin.  In the next theorem we consider
the degenerate case $\sigma^2=0$ in one dimension.
\medskip

\noindent{\bf Theorem 3} {\it Let $K(x,y)=A(x-y)$ be a
translation-invariant kernel in $\bbR^1$ and $m(\lambda )
=\vert\lambda\vert^\alpha\varphi (\lambda )$, where $\varphi
(\lambda )$ is a slowly varying function at the origin and
$0<\alpha <1$.  Then for any Schwartz function $f$ :
$E S_{f_L}=LA(0)\int f(x)dx$, Var $S_{f_L}=L^{1-\alpha}
\varphi (L^{-1})\int\vert\hat f(k)\vert^2\vert k
\vert^\alpha dk (1+o(1))$, and 
$$\frac{S_{f_L}-ES_{f_L}}
{(L^{1-\alpha}\varphi (L^{-1}))^{\tfrac{1}{2}}}$$ 
converges
in distribution to $N(0,\int\vert\hat f(k)\vert^2\vert
k\vert^\alpha dk$.}
\medskip

\noindent{\bf Remark 5}  We recall that $\varphi (\lambda) \geq 0$
is slowly varying at the origin if $\lim_{\lambda\rightarrow
0}\tfrac{\varphi (a\lambda )}{\varphi (\lambda )}=1$ for
any $a\neq 0$ (see [Se]).
\medskip

\noindent{\bf Remark 6}  The result of Theorem 3 can be
interpreted as the convergence in distribution of the
random signed measure
$$\biggl (L^{1-\alpha }\varphi (L^{-1})\biggr )
^{-\tfrac{1}{2}}\left (
\sum\delta \left (x-\frac{x_i}{L}\right )-A(0)Ldx\right )$$
to the self-similar (also called automodel in the Russian literature) 
generalized Gaussian random
process with the spectral density $\vert k\vert^\alpha, \ 0<\alpha <1.$
Self-similarity  means
that the distribution of the process is invariant under the
action of the renorm-group $\xi (x)\rightarrow\xi (ax)a^\gamma$,
$\gamma =\tfrac{1+\alpha}{2}$ . The self-similar generalized Gaussian random 
process corresponding to $\alpha =0$ is exactly the white noise (see Remark 3
above). It was proven by Dobrushin that the only self-similar random processes
in $\bbR^1$ are the ones with the spectral density
$\vert k\vert^\alpha, \ 0\leq\alpha \leq 1 $. 
A self-similar generalized random process with the spectral density $|k| \ $
appeared in the Spohn's results ([Sp]) discussed above after the formulas 
$ (4), \ (5)$ (see also [Jo1], [B], [So2]).
For additional information 
on self-similar random processes we refer the reader to [Dob1], [Dob2],[S].
\medskip

\noindent{\bf Example}  Let $\hat A$ be the indicator of 
$\sqcup_{n\geq 1}[n,n+n^{-\beta}],\ \beta >1$.  Then
$m(\lambda )=\const\cdot\vert\lambda\vert^{1-\tfrac{1}
{\beta}}(1+o(1))$. On the other hand if the length
$l_n$  of the $n-$th interval $ [n,n+l_n]$ decays 
sufficiently fast, say $ 0 \leq l_{n+1} \leq  l_n^{1+\epsilon}, \ \ 
\epsilon >0 \ $ than $m(\lambda)$  is not regularly varying
at the origin.

Finally we consider the case when $\hat A$ is the indicator
of a union of $1\leq\ell <\infty$ disjoint intervals.
It is straightforward to see that then $m(\lambda )=\ell
\vert\lambda\vert$ near the origin.
\medskip

\noindent{\bf Theorem 4} {\it Let $\hat A$ be the indicator
of $I$, $I=\sqcup^\ell_{i=1}[a_i,b_i],\ a_1<b_1<a_2<b_2<
\dots <a_\ell <b_\ell$. Then for any Schwartz function $f \ \ $
$\sum^\infty_{i=-\infty}f\left (\tfrac{x_i}{L}\right )
-A(0)\cdot L\int^\infty_{-\infty}f(x)dx$ converges in
distribution to $N(0,\ell\cdot\int^\infty_{-\infty}\vert
\hat f(k)\vert^2\vert k\vert dk)$.}
\medskip

The proofs of Theorems 1--3 will be given in the next three
sections.  The proof of Theorem 4 is the same, modulo
trivial alterations, as the one given for the sine kernel
in [So2].

The author would like to thank P.Deift and Ya. Sinai for useful discussions,
T. Shirai 
for making available the preprint [ST] prior to its publication and the referee
for useful comments.

\section{Proof of Theorem 1}

We are going to prove Theorem 1 by the method of moments.
Let us denote by $C_n(S_f)$ the $n^{th}$ cumulant of $S_f$.
We remind the reader that for a random variable $\eta $ with all finite moments
, the cumulants $ C_n (\eta), \ n=1,2,\ldots$ are defined through the Taylor
coefficients of the logarithm of the characteristic function :
$$ \log E( \exp(it\eta)) = \sum_{n=1}^{\infty} C_n(\eta) (it)^n/n! .$$
We show that the $n^{th}$ cumulant of the normalized linear
statistics $\tfrac{S_{f_L}-ES_{f_L}}{({\rm Var}\ S_{f_L}
)^{\tfrac{1}{2}}}$ converges to zero as $L\rightarrow
\infty$ for sufficiently large $n$ ($n>\max (2\delta ,2))$.
The Lemma 3 from the Appendix then asserts that all
cumulants of $\tfrac{S_{f_L}-ES_{f_L}}{({\rm Var}\ S_{f_L}
)^{\tfrac{1}{2}}}$ converge to the cumulants of the standard
normal distribution, which implies the weak convergence.

We recall the lemma established in [So2] (see formula 
(2.7)).
\medskip

\noindent{\bf Lemma 1}
\begin{equation}
\begin{split}
&C_n(S_f)=\sum^n_{m=1}\sum_{\overset{(n_1,\dots ,n_m):
n_1+\dots +n_m=n.}{n_1\geq 1,\ i=1,\dots ,m}}\frac
{(-1)^{m-1}}{m}\frac{n!}{n_1!\dots n_m!}
\int f^{n_1}(x_1)\\
&K(x_1,x_2)f^{n_2}(x_2)K(x_2,x_3)\dots
f^{n_m}(x_m)K(x_m,x_1)d\mu (x_1)\dots d\mu (x_m)
\end{split}
\end{equation}
\medskip

Using Lemma 1 we will be able to estimate the cumulants
of $S_{f_L}$.  We claim the following result to be true.
\medskip

\noindent{\bf Lemma 2} {\it Under the assumptions of 
Theorem 1
\begin{equation}
C_n(S_{f_L})=O(({\rm Var}_LS_{ f_L})^{\delta +\epsilon}),\ n\geq 1, 
\end{equation}
 where $ \epsilon$ is arbitrarily small.}
\medskip

\noindent{\bf Proof of Lemma 2} It follows from (14) that
$C_n(S_{f_L})$ is a linear combination of
\begin{equation*}
\begin{split}
&\int f^{n_1}_L(x_1)K_L(x_1,x_2)f^{n_2}_L(x_2)K_L(x_2,x_3)\dots
f^{n_m}_L(x_m)K_L(x_m,x_1)\\
&d\mu (x_1)\dots d\mu (x_m)=\tr (f^{n_1}_LK_Lf^{n_2}_LK_L\dots
f^{n_m}_LK_L),
\end{split}
\end{equation*}

where $n_i\geq 1,\ i=1,\dots ,m,\ m\geq 1.$

We claim that each term is $O(({\rm Var}_LS_{ f_L})^{\delta +\epsilon}).$
Indeed, if $m=1$, then \\
$\vert\tr f^n_LK_L\vert=\vert\int
f^n_L(x)K_L(x,x)d\mu (x)\vert \leq \Vert f_L\Vert^{n-1}_\infty
\int\vert f_L(x)\vert K_L(x,x)d\mu (x)=\Vert f_L
\Vert^{n-1}_\infty ES_{\vert f\vert_L} =
O(( {\rm Var}_LS_{ f_L})^{\delta +\epsilon}).$

If $m>1$, 
represent $\tr (f^{n_1}_LK_Lf^{n_2}_L
K_L\dots f^{n_m}_LK_L)$ as a linear combination of $\tr 
(f^{n_1}_{\pm ,L}K_Lf^{n_2}_{\pm ,L}K_L
\dots f^{n_m}_{\pm ,L}
K_L)$, where we use the notations $f_+=\max (f,0)$, $f_-=\max (-f 
,0)$. Let us fix the
choice of $\pm$ in each of the factors.  Using the cyclicity
of the trace and the inequality $\vert\tr (AB)\vert\leq (
\tr (AA^*))^{\tfrac{1}{2}}(\tr (BB^*))^{\tfrac{1}{2}}$
for the Hilbert-Schmidt operators ([RS], section VI.6),
we obtain
\begin{equation}
\begin{split}
&\vert\tr (f^{n_1}_{\pm ,L}K_Lf^{n_2}_{\pm ,L}K_L\dots 
f^{n_m}_{\pm ,L}K_L)\vert =
\vert\tr (f^{\frac{n_1}{2}}_{\pm ,L}K_Lf^{n_2}_{\pm ,L}K_L
\dots f^{n_m}_{\pm ,L}K_Lf^{\frac{n_1}{2}}_{\pm ,L})\vert
\leq\\
&[\tr ((f^{\frac{n_1}{2}}_{\pm ,L}K_Lf^{\frac{n_2}{2}
}_{\pm ,L})(f^{\frac{n_1}{2}}_{\pm ,L}K_Lf^{\frac{n_2}{2}
}_{\pm ,L})^*)]^{\frac{1}{2}}[\tr ((f^{\frac{n_2}{2}
}_{\pm ,L}K_Lf^{n_3}_{\pm ,L}K_L\dots f^{n_m}_{\pm ,L}K_L
f^{\frac{n_1}{2}}_{\pm ,L})\\
&(f^{\frac{n_2}{2}}_{\pm ,L}K_Lf^{n_3}_{\pm ,L}K_L\dots
f^{n_m}_{\pm ,L}K_Lf^{\frac{n_1}{2}}_{\pm ,L})^*)]^{\frac
{1}{2}}.
\end{split}
\end{equation}
The first factor at the r.h.s. of (16) is equal (again
by the cyclicity of the trace) to $[\tr (f^{n_1}_{\pm ,L}
K_Lf^{n_2}_{\pm ,L}K_L)]^{\tfrac{1}{2}}$ (in particular we
note that $\tr (g_1Kg_2K)\geq 0$ for non-negative $g_1,\
g_2$).

Since
\begin{equation*}
\begin{split}
&\tr ((f^{n_1}_{\pm ,L}+f^{n_2}_{\pm ,L})^2K)-\tr
((f^{n_1}_{\pm ,L}+f^{n_2}_{\pm ,L})K(f^{n_1}_{\pm ,L}
+f^{n_2}_{\pm ,L})K)=\\
&{\rm Var}\ (S_{f^{n_1}_{\pm ,L}+f^{n_2}_{\pm ,L}})\geq 0,
\end{split}
\end{equation*}
we have
\begin{equation}
\begin{split}
&0\leq\tr (f^{n_1}_{\pm ,L}K_Lf^{n_2}_{\pm ,L}K_L)
\leq\frac{1}{2}\biggl (\tr ((f^{n_1}_{\pm ,L}+f^{n_2}_{\pm ,L})^2K_L)-
\tr (f^{n_1}_{\pm ,L}K_Lf^{n_1}_{\pm ,L}K_L)-\\
&\tr (f^{n_2}_{\pm ,L}K_Lf^{n_2}_{\pm ,L}K_L)\biggr )\leq\frac{1}{2}
\tr ((f^{n_1}_{\pm ,L}+f^{n_2}_{\pm ,L})^2K_L))=O(\tr (\vert f_L
\vert K_L)) o(({\rm Var}_L S_{f_L})^{\epsilon})\\
&=O(({\rm Var}_L S_{f_L})^{\delta +\epsilon}).
\end{split}
\end{equation}

As for the second term in (16), one can rewrite $\tr\biggl ((
f^{\tfrac{n_2}{2}}_{\pm ,L}K_Lf^{n_3}_{\pm ,L}K_L\dots
\break f^{n_m}_{\pm ,L}K_Lf^{\tfrac{n_1}{2}}_{\pm ,L})
(f^{\tfrac{n_2}{2}}_{\pm ,L}K_Lf^{n_3}_{\pm ,L}K_L\dots
f^{n_m}_{\pm ,L}K_Lf^{\tfrac{n_1}{2}}_{\pm ,L})^*\biggr )$ as
\begin{equation}
\tr\biggl (f^{\frac{n_2}{2}}_{\pm ,L}K_Lf^{n_3}_{\pm ,L}K_L\dots
f^{n_m}_{\pm ,L}K_Lf^{n_1}_{\pm ,L}K_Lf^{n_m}_{\pm ,L}\dots
K_Lf^{n_3}_{\pm ,L}K_Lf^{\frac{n_2}{2}}_{\pm ,L}\biggr )=\tr (CDD^*),
\end{equation}
where $C=f^{\tfrac{n_3}{2}}_{\pm}K_Lf^{n_2}_{\pm ,L}K_L
f^{\tfrac{n_3}{2}}_{\pm ,L}$, $D=f^{\tfrac{n_3}{2}}_{\pm ,L}
K_Lf^{n_4}_{\pm ,L}K_L\dots f^{n_m}_{\pm ,L}K_L f^{\tfrac{n_1}{2}}_{\pm ,L}$.
  Note that $C\geq 0$ and
$\tr (C)=\tr (f^{n_3}_{\pm ,L}K_Lf^{n_2}_{\pm ,L}K)=O
( ({\rm Var}_L S_{f_L})^{\delta +\epsilon})$ by  arguments similar to
(17).  Using $\vert\tr (CDD^*)\vert\leq\tr (C)\cdot\Vert
DD^*\Vert=\tr (C)\cdot\Vert D\Vert^2$ ([RS], Section
VI.6) and $\Vert D\Vert\leq\Vert K\Vert^m\cdot\Vert f_L
\Vert^\aleph_\infty$, where $\aleph =(\sum^m_{i=1}n_1)-n_2$,
we conclude that (18) is $O( ({\rm Var}_L S_{f_L})^{\delta +\epsilon})$
Together with (16) and (17) this concludes the proof of
the lemma.\hfill$\dashv$

Let us now apply Lemma 2 to estimate the cumulants of the
normalized linear statistics.  We have 
$$C_1\left (\frac{S_{f_L}-ES_{f_L}}{\sqrt{\var_L S_{f_L}}}
\right )=0,$$
$$C_2\left (\frac{S_{f_L}-ES_{f_L}}{\sqrt{\var_L\ S_{f_L}}}
\right )=1,$$
and, for $n>2$,
\begin{equation}
C_n\left (\frac{S_{f_L}-ES_{f_L}}{\sqrt{\var_L\ S_{f_L}}}
\right )=\frac{C_n(S_{f_L})}{(\var_L\ S_{f_L})^{\frac{n}{2}}}
=O\left (\frac{E\left (S_{\vert f\vert_L}\right )}
{\left (\var_L\ S_{f_L}\right )^{\frac{n}{2}}}\right )
\end{equation}
It follows from the Lemma 2 and (19)  that $C_n\left (
\tfrac{S_{f_L}-ES_{f_L}}{\sqrt{\var\ S_{f_L}}}\right )$ 
goes to zero if $n>2\delta$.

Lemma 3 from the Appendix then implies that
all cumulants of the normalized linear statistics converge
to the cumulants of the standard normal random variable,
and weak convergence of the distributions follows.

Theorem 1 is proven.\hfill$\Box$

\section{Proof of Theorem 2}

Let $(E,d\mu )$ be $(\bbR^d,dx)$ and $T_Lx=\tfrac{x}{L}$.
Consider a real-valued function $f\in L^1(\bbR^d)\cap L^2
(\bbR^d)$.  The mathematical expectations of $S_{f_L}$
is equal to
\begin{equation*}
\begin{split}
&ES_{f_L}=\int_{\bbR^d}f(x /L )K_L(x,x)dx=
\int_{\bbR^d}f(x/L )A_L(0)dx\\
&+\int_{\bbR^d}f(x/L )R_L(x,x)dx=A_L(0)\cdot
L^d\int f(x)dx+\int f(x/L )R_L(x,x)dx.
\end{split}
\end{equation*}
By (11) the absolute value of the second integral is bounded
by the sum of the integrals
\begin{equation*}
\begin{split}
&\int_{\bbR^d_+}\left \vert f(\pm x_1/L,\dots
,\pm x_d/L )\right \vert Q\left (2x_1,\dots ,
2x_d\right )dx\\
&\leq\left (\int_{\bbR^d_+}f^2(\pm x_1/L,
\dots ,\pm x_d/L)dx\right )^{\frac{1}{2}}
\left (\int_{\bbR^d_+}Q^2(2x)dx\right )^{\frac{1}{2}}
=O\left (L^{\frac{d}{2}}\right )
\end{split}
\end{equation*}
Therefore,
\begin{equation}
ES_{f_L}=A_L(0)\cdot L^d\cdot\int_{\bbR^d}f(x)dx+O
(L^{\frac{d}{2}}).
\end{equation}
The variance of $S_{f_L}$ is given by
\begin{equation}
\begin{split}
&\var\ S_{f_L}=\int f^2(x/L)K_L(x,x)dx
-\int f (x/L )f(y/L )
\vert K_L(x,y)\vert^2dxdy=\\
&A_L(0)L^d\int f^2(x)dx-\int f(x/L )
f(y/L )\vert A_L(x-y)\vert^2dxdy+ r(L),
\end{split}
\end{equation}
where
\begin{equation*}
\begin{split}
&r(L)=\int f^2(x/L)R_L(x,x)dx-
2\int f(x/L)f(y/L )
A_L(x-y)R_L(y,x)dxdy-\\
&\int f(x/L )f(y/L )
\vert R_L(x,y)\vert^2dxdy=r_1(L)+r_2(L)+r_3(L).
\end{split}
\end{equation*}
It follows from the assumptions of the theorem that the second term at the 
r.h.s. of (21) is equal to
\begin{equation*}
\begin{split}
&L^d\int\vert\hat f(k)\vert^2\widehat{\vert A_L\vert^2}
(k/L )dk=L^d\widehat{\vert A_L\vert^2}
(0)\int\vert\hat f(k)\vert^2dk\cdot (1+o(1))=\\
&L^d\int\vert A_L(x)\vert^2dx\int f^2(x)dx(1+o(1)).
\end{split}
\end{equation*}
Indeed,
\begin{equation*}
\begin{split}
&\vert\int\vert\hat f(k)\vert^2 ( \widehat{\vert A_L\vert^2}
(k/L ) -\widehat{\vert A_L\vert^2}(0))  dk \vert \leq \\
& \vert \int_{|k|> \kappa_L} \vert +\vert \int_{|k|\leq \kappa_L} \vert.
\end{split}
\end{equation*}
Since
\begin{equation*}
\begin{split}
&\vert \widehat{\vert A_L\vert^2}(t)\vert = \vert \int \hat{A_L}(k) 
\hat{A_L}(k-t)dk \vert \\
& \leq \int \hat{A_L}(k) dk =A_L(0) \leq const
\end{split}
\end{equation*}
we note that the first integral is bounded from above by 
$$ const \int_{|k|>(\kappa_L)^{1/2}} |\hat{f}(k)|^2 dk \ \ \to 0 \ \ {\rm as}
\ \ L \to \infty .$$
To deal with the second integral we estimate from above 
\begin{equation*}
\begin{split}
&\vert \widehat{\vert A_L\vert^2}
(k/L ) -\widehat{\vert A_L\vert^2}(0)\vert  \leq \\
& \vert \int |A_L|^2(t) (\exp(2\pi i t k/L) -1) dt \vert \\
& \vert \int_{|t|\geq L/\kappa_L}  + \int_{ |t| < L/\kappa_L} \vert \leq \\
& \int_{|t|\geq L/\kappa_L}  |A_L|^2(t) dt + O(1/\sqrt(\kappa_L))= \\
& o(1) + O(1/\sqrt(\kappa_L))=o(1).
\end{split}
\end{equation*}
Therefore
\begin{equation}
\begin{split}
&\var\ S_{f_L}=\left (A_L(0)-\int\vert A_L(x)\vert^2dx\right ) 
L^d\int f^2(x)dx+o(L^d)+r(L)=\\
&\sigma^2L^d\int f^2(x)dx+o(L^d)+r(L).
\end{split}
\end{equation}
We claim that
\begin{equation}
r(L)=o(L^d).
\end{equation}
Consider first $r_1(L)$.  By (11) it is bounded by the
integrals
$$\int_{\bbR^d_+}f^2 (\pm x_1/L,\dots
,\pm x_d/L )Q(2x)dx$$
All of these integrals are estimated in the same way.
For example,
\begin{equation*}
\begin{split}
&\int_{\bbR^d_+}f^2 (x/L )Q(2x)dx=L^d
\int f^2(x)Q(2Lx)dx=\\
&L^d\int f^2(x)Q(2Lx)\x_{\{Q(2Lx)>\frac{1}{\sqrt{L}}\}}
dx+L^d\int f^2(x)Q(2Lx)\x_{\{Q(2Lx)\leq\frac{1}{\sqrt{L}}
\}}\\
&\leq L^d\Vert Q\Vert_\infty\int f^2(x)\x_{\{Q(2Lx)>\frac
{1}{\sqrt{L}}\}}dx+L^{d-\frac{1}{2}}\int f^2 (x)dx=o(L^d),
\end{split}
\end{equation*}
since
$$\ell (x:Q(2Lx)>\frac{1}{\sqrt{L}})\underset
{L\rightarrow\infty}{\longrightarrow}0.$$
To estimate $r_3(L)$ we need to estimate the integrals
of the form
\begin{equation}
\int_{\bbR^d_+}\vert f(x/L )
\vert \ \vert f(y/L )
\vert  Q^2(x+y)dxdy=L^d\int_{\bbR^d_+}g (
z/L)Q^2(z)dz,
\end{equation}
where $g(z)=\int\vert f(x)\vert\ \vert f(z-x)\vert
\x_{\bbR^d_+}(x)\x_{\bbR^d_+}(z-x)dx$.

Since $g(z)$ is bounded, continuous, and zero at the origin,
we have
$$(24)=L^dg(0)\int Q^2(z)dz (1+o(1))=o(L^d).$$
Finally,
\begin{equation*}
\begin{split}
&\vert r_2(L)\vert=\vert \ \int f (x/L
)f (y/L)A_L(x-y)R(y,x)dxdy
\vert\\
&\leq\left [\int \ \vert f(x/L )
\vert \ \vert f(y/L )
\vert\ \vert A_L(x-y)\vert^2dxdy\right ]^{\frac{1}{2}}\\
&\left [\int \ \vert f(x/L )
\vert\ \vert f(y/L )\vert\ 
\vert R_L(y,x)\vert^2dxdy\right ]^{\frac{1}{2}}=O\left (
L^{\frac{d}{2}}\right )\ o\left (L^{\frac{d}{2}}\right )
=o\left (L^d\right ).
\end{split}
\end{equation*}
Combining the above estimates, we prove (23), which implies
\begin{equation}
\var\ S_{f_L}=\sigma^2L^d\int f^2(x)dx(1+o(1)).
\end{equation}
If $f$ is bounded, the Central Limit Theorem then follows
from Theorem 1 (compactness of the support of $f$ is not
needed since $f\in L^1(\bbR^d)\cap L^\infty (\bbR^d)$
guarantees that all moments of $S_{f_L}$ are finite).  The
proof in the case of the unbounded $f$ follows by a rather
standard approximation argument.  We choose $N>0$ to be
sufficiently large and consider a truncated function
$$\tilde f(x)=\begin{cases}f(x),&\text{if }\vert f(x)
\vert\leq N\\N,&\text{if }f(x)>N\\-N&\text{if }
f(x)<-N.\end{cases}$$
Observe that
$$E\left (\frac{S_{f_L}-ES_{f_L}}{\sigma L^{\frac{d}{2}}
}-\frac{S_{\tilde f_L}-ES_{\tilde f_L}}{\sigma L^{\frac
{d}{2}}}\right )^2=\frac{\var\ S_{(f-\tilde f)_L}}
{\sigma^2L^d}=\frac{\int_{\vert x\vert\geq N}f^2(x)dx}
{\sigma^2}+o(1)$$
can be made arbitrarily small by choosing $N$ and $L$ 
sufficiently large.

Since
$$\frac{S_{\tilde f_L}-ES_{\tilde f_L}}{\sigma L^{\frac{d}
{2}}}\overset{w}{\underset{L\rightarrow\infty}
{\longrightarrow}}N\left (0,\int_{\vert x\vert\leq N}f^2(x)dx
\right )$$
and
$$\lim_{N\rightarrow\infty}\int_{\vert x\vert\leq N}f^2
(x)dx=\int f^2(x)dx,$$
the  result follows.

Theorem 2 is proven.\hfill$\Box$

\section{Proof of Theorem 3}
We now turn to the proof of Theorem 3.
It is enough to establish that
\begin{equation}
\var\ S_{f_L}=L^{1-\alpha}\varphi (L^{-1})\int\vert\hat f
(k)\vert^2\ \vert k\vert^\alpha dk(1+o(1)).
\end{equation}
The result then will follow from Theorem 1.  We have (see (12))
\begin{equation}
\begin{split}
&\var\ S_{f_L}=\int\vert\hat f (L\lambda )\vert^2L^2m
(\lambda )d\lambda =L\int\vert\hat f(k)\vert^2m(kL^{-1})
dk=\\
&L\int\vert\hat f(k)\vert^2\ \vert k\vert^\alpha L^{-\alpha}
\varphi (kL^{-1})dk=L^{1-\alpha}\varphi (L^{-1})\\
&\int\vert
\hat f(k)\vert^2\ \vert k\vert^\alpha\frac{\varphi 
(kL^{-1})}{\varphi (L^{-1})}dk
\end{split}
\end{equation}
It was proven by Karamata ([K1], [K2]) that any slowly
varying function at the origin can be represented in some
interval $(0,b]$ as
\begin{equation}
\varphi (x)=\exp\left \{\eta (x)+\int^{x^{-1}}_{b^{-1}}\frac
{\epsilon (t)}{t}dt\right \},
\end{equation}
where $\eta$ is a bounded measurable function on $(0,b]$,
such that $\eta (x)\rightarrow c$ as $x\rightarrow 0$
($\vert c\vert <\infty$), and $\epsilon (x)$ is a continuous
function on $(0,b]$ such that $\epsilon (x)\rightarrow 0$
as $x\rightarrow 0$.  (for a modern day reference we refer the reader to [Se],
Theorem 1.2; of course a similar representation
holds for $\varphi$ also on some interval $[b', 0)$ of the
negative semi-axis).  In particular
\begin{equation}
\frac{\varphi\left (\frac{k}{L}\right )}
{\varphi\left (\frac{1}{L}\right )}
\underset{L\rightarrow\infty}{\longrightarrow}1
\end{equation}
uniformly in $k$ on compact subsets of $\bbR^1\setminus\{0\}$,
and the following estimates hold uniformly in $k$ for 
sufficiently large $L$
\begin{equation}
\const_1k^{-n}\leq\varphi (kL^{-1})/\varphi (L^{-1})\leq
\const_2k^n,\text{ for }1\leq k\leq L,
\end{equation}
\begin{equation}
\const_3k^{-\frac{1}{2}}\leq\varphi (kL^{-1})/\varphi
(L^{-1})\leq\const_4k^{\frac{1}{2}}\text{ for }0<k\leq 1,
\end{equation}
where $\const_i,\ i=1,\dots 4,\ n>0$ are some constants.

The estimates (28)-(31) imply
$$\int^L_{-L}\vert\hat f(k)\vert^2\ \vert k\vert^2\frac
{\varphi (kL^{-1})}{\varphi (L^{-1})}dk
\underset{L\rightarrow\infty}{\longrightarrow}
\int^\infty_{-\infty}\vert\hat f(k)\vert^2\ \vert k
\vert^\alpha dk.$$
From the other side, the integral over $\vert k\vert\geq 
L$ is o(1) since $f$ is a Schwartz function and $m$ is
bounded.

Theorem 3 is proven.\hfill$\Box$

\noindent{\bf Remark 7} 
We learned very recently that similar results to our Theorem 2  
have been independently
obtained (in the discrete case) by Tomoyuki Shirai and Yoichiro 
Takahashi in the
preprint [ST].

\section*{Appendix}

For the convenience of the reader we give here the proof
of a rather standard fact.
\medskip

\noindent{\bf Lemma 3} {\it Let $\{\eta_L\}$ be a family
of random variables such that $c_1(\eta_L)=0,\ c_2
(\eta_L)=1$ and $c_n(\eta_L)$ converges to zero as $L
\rightarrow\infty$ for all $n\geq N$, where $N<\infty$.
Then $\lim_{L\rightarrow\infty}c_n(\eta_L)=0$ for all
$n>2$ and $\eta_L$ converges in distribution to $N(0,1)$.}
\medskip

\noindent{\bf Proof}  Denote $d_L=\max (\vert c_j(\eta_L)
\vert^{\tfrac{1}{j}},\ 1\leq j\leq N-1)$.  It is clear that
$d_L\geq 1$.  Consider the random variable
$$\tilde\eta_L=\eta_L/d_L.$$
Since $c_n(\tilde\eta_L)=c_n(\eta_L)/d^n_L$ we have
$\vert c_n(\tilde\eta_L)\vert\leq 1$ for all $n$ and
$c_n(\tilde\eta_L)\rightarrow 0$ for $n\geq N$.
Consider $(N-1)$-dimensional vector $(c_1(\tilde\eta_L),
\dots ,c_{N-1}(\tilde\eta_L))$.  Let $(c_1,c_2,\dots ,
c_{N-1})$ be a limit point.  The Marcinkiewicz theorem (see e.g. [L])
states that if all but a finite number of cumulants of a random variable
are non-zero then  the random variable must either have a Gaussian distribution
or be a constant. In both cases 
we have $c_j=0$ for $j>2$.  Therefore $d_L=(c_2(\eta_L
))^{\tfrac{1}{2}}=1$ for sufficiently large $L$ and 
$c_n(\eta_L)\underset{L\rightarrow\infty}{\longrightarrow}
0$ for $n>2$.  Convergence of the cumulants of $\eta_L$
to the cumulants of $N(0,1)$ is equivalent to the 
convergence of the moments which in turn implies convergence
in distribution.

\def\am{{\it Ann. of Math.} }
\def\ap{{\it Ann. Probab.} }
\def\temf{{\it Teor. Mat. Fiz.} }
\def\jmp{{\it J. Math. Phys.} }
\def\cmp{{\it Commun. Math. Phys.} }
\def\jsp{{\it J. Stat. Phys.} }
\def\prl{{\it Phys. Rev. Lett.} }
\def\arma{{\it Arch. Rational Mech. Anal.} }

\end{document}